\documentclass[font=11pt]{amsart}
\usepackage{amstext,amsmath,amssymb,amsthm}
\usepackage{latexsym}
\usepackage{exscale}
\usepackage{tikz}
\usepackage{xcolor}
\usepackage[margin=1.25in]{geometry}
\usepackage{hyperref}
\usepackage[capitalize]{cleveref}
\usepackage{enumerate}

\numberwithin{equation}{section}

\renewcommand\r{\rangle}
\renewcommand\l{\langle}
\newcommand\dsize{\displaystyle}

\newcommand\dist{\operatorname{dist}}

\renewcommand\epsilon{\varepsilon}

\newcommand\R{\mathbb{R}}
\newcommand\C{\mathbb{C}}

\newcommand\Z{\mathbb{Z}}
\newcommand\T{\mathbb{T}}
\newcommand\N{\mathbb{N}}

\newcommand\B{\mathcal{B}}

\newcommand\calV{\mathcal{V}}

\newcommand\bA{\boldsymbol{A}}
\newcommand\bB{\boldsymbol{B}}
\newcommand\bE{\boldsymbol{E}}
\newcommand\bF{\boldsymbol{F}}
\newcommand\bI{\boldsymbol{I}}
\newcommand\bU{\boldsymbol{U}}
\newcommand{\bV}{\boldsymbol{V}}
\newcommand\G{\boldsymbol{\Gamma}}

\newcommand{\andspace}{\quad\text{and}\quad}
\newcommand{\wherespace}{\quad\text{where}\quad}
\newcommand{\sep}{{\rm sep}}
\newcommand\sign{{\rm sign}}

\newtheorem{Thm}{Theorem}[section]

\newtheorem{Cor}[Thm]{Corollary}
\newtheorem{Prop}[Thm]{Proposition}
\theoremstyle{remark}
\newtheorem{Rem}{Remark}[section]

\begin{document}

\title[]{Concerning the stability of exponential systems and Fourier matrices  }

\author{Oleg Asipchuk}
\address{Oleg Asipchuk, Florida International University, Miami (FL)}
\email{aasip001@fiu.edu }

\author{Laura De Carli}
\address{Laura De Carli, Florida International University, Miami (FL)}

\email{decarlil@fiu.edu }

\author{Weilin Li}
\address{Weilin Li. City University of New York, City College (NY)}
\email{wli6@ccny.cuny.edu}

\begin{abstract}
	Fourier matrices naturally appear in many applications and their stability is closely tied to performance guarantees of algorithms. The starting point of this article is a result that characterizes properties of an exponential system on a union of cubes in $\R^d$ in terms of a general class of Fourier matrices and their extreme singular values. This relationship is flexible in the sense that it holds for any dimension $d$, for many types of exponential systems (Riesz bases, Riesz sequences, or frames)  and for Fourier matrices with an arbitrary number of rows and columns. From there, we prove new stability results for Fourier matrices by exploiting this connection and using powerful stability theorems for exponential systems. This paper provides a systematic exploration of this connection and suggests some natural open questions.
\end{abstract} 

\maketitle

\section{Introduction}

The main purpose of this paper is to exploit a connection between exponential systems on  unions of hyper-cubes of  $\R^d$ and their associated Fourier matrices. From this relationship, we transfer sophisticated results from one subject to the other.

We start with some preliminary notation, and further background is provided in \cref{sec:prelim}. Throughout, we let $d\ge 1$ be an integer and $Q:=Q_d= [0, 1)^d$. Let  $\Delta:=\{\vec \delta_j\}_{j=1}^L$ be a subset of $\mathbb{T}^d:=(\R/\Z)^d$ and $P:=\{\vec p_k\}_{k=1}^N \subset \Z^d$ consists of distinct elements. We consider the family of exponential functions, 
\begin{equation}\label{e-BDelta}
	\B
	:=\B (\Delta)
	=: \bigcup_{j=1}^L \big\{e^{ 2\pi i (n+\vec \delta_j)\cdot x} \big\}_{n\in\Z^d},
\end{equation}
defined on the following union of non-overlapping cubes,
\begin{equation}\label{eq:TP}
	T=T(P)=:  \bigcup_{k=1}^N \, (Q  +\vec p_k).
\end{equation} 

A classical question in harmonic analysis is to determine properties of the exponential set $\B(\Delta)$  on the domain $T(P)$. The starting point of this paper is a result which shows that this question can be equivalently recast purely in terms of a single matrix. We associate  the exponential set $\B(\Delta)$ on the set $T(P)$ with the $L\times N$ complex matrix
\begin{equation}\label{def-Gamma} 
	\G
	:=\G(\Delta,P)
	:= \Big[ e^{2\pi i \delta_j\cdot p_k} \Big]_{j=1,\dots,L \atop{  k=1,\dots, N}}\ . 
\end{equation}
Since we have placed no assumptions on $\Delta\subset\T^d$ and $P\subset\Z$, this matrix could be highly irregular. Under special circumstances, it reduces to a Vandermonde, non-harmonic Fourier, nonuniform Fourier, or partial discrete Fourier transform matrix (these definitions are recalled in \cref{2.3}).

The following theorem relates the frame property of $\B(\Delta)$ for $L^2(T(P))$ and extreme singular values of $\G(\Delta,P)$. It is a generalization of the basis/square case, where it is known (see \cite[Thm 1.1]{DC},\cite{BK, GL, K2, AAC}) that $\B(\Delta)$ is a Riesz basis for $L^2(T(P))$ if and only if $L=N$ and the square matrix $\G(\Delta,P)$ is invertible.  Some results  related to non-square matrices are also in \cite{NOU2}.

\begin{Thm}\label{T1-frame}
	Let $L=|\Delta|$ and $N=|P|$. 
	\begin{enumerate}[\rm (a)]
		
		\item 
		Assume $L  \ge N$. The set 
		$\B (\Delta)$ is a frame for $L^2(T(P))$  if and only if  the associated matrix $\G(\Delta,P)$ has rank $N$.  Equivalently, $\B(\Delta)$ is a frame for $L^2(T(P))$ if and only if the rows  of $\G(\Delta,P)$ form a frame for $\C^N$. If any of these equivalent statements hold, the optimal frame constants of $\B(\Delta)$ for $L^2(T(P))$ are $\sigma_1^2(\G(\Delta,P) $ and $\sigma_N^2(\G(\Delta,P) $.
		\item 
		Assume $L\leq N$. The set $\B (\Delta)$ is a Riesz sequence in $L^2(T(P))$ if and only if  the associated matrix $\G(\Delta,P)$ has rank $L$.  Equivalently, $\B(\Delta)$ is a Riesz sequence in $L^2(T(P))$ if and only if the rows of $\G(\Delta,P)$ are linearly independent. If any of these equivalent statements hold, the optimal Riesz sequence constants of $\B(\Delta)$ in $L^2(T(P))$ are $\sigma_1^2(\G(\Delta,P) $ and $\sigma_L^2(\G(\Delta,P) $.
	\end{enumerate}	
\end{Thm}

The proof of \cref{T1-frame} can be found in \cref{sec:proofframe}. For part (a), the condition $L\leq N$ is necessary for the conclusion to hold and $\G$ is a tall rectangular matrix. Similarly, for part (b), the condition $L\geq N$ is necessary and $\G$ is a wide rectangular matrix. This theorem implies that when $L=N$, the system $\B(\Delta)$ is a Riesz sequence in $L^2(T(P))$ if and only if it is a frame for $L^2(T(P))$.

It  also shows that any result for exponential systems on unions of cubes has a corresponding statement for the extreme singular values of its associated matrix, and vice versa. We  cannot flush out all ramifications of this result in a single paper, so we focus on a more specific question: the stability of exponential systems has been well studied (see e.g., the classical Kadec theorem and its variations), but what are the corresponding statements for  the matrices associated to these systems? 

Our investigation of this question revealed that the corresponding statements describe how the extreme singular values of a Fourier matrix $\bF(\Omega,U)$ change when its nodes $U$ and/or frequencies $\Omega$ are perturbed by a small amount to $U'$ and/or $\Omega'$. We refer to them as {\it node and frequency stability}, respectively. 

Just to give the reader a flavor of  the results that will be presented in this paper, we state a special case of Theorem \ref{T3-freqstability} in Section 3, where    we prove     stability results for Fourier matrices associated to exponential systems  on rectangles of $\R^d$.
  
\begin{Thm}\label{T-Kad-Matr} 
	Let $\bV$  be the $N\times N$  Vandermonde matrix with nodes $ z_k= e^{2\pi i  (k+\epsilon_k)/N} $, where $\ell:=\max_k |\epsilon_k| <\frac 14$. Then $\bV$  is invertible and its extreme singular values satisfy the inequalities $\sigma_N\ge  N(\cos(\pi \ell)-\sin( \pi \ell))$ and	$\sigma_1\leq N(2-\cos(\pi \ell)+\sin( \pi \ell)) $.
\end{Thm}

Our paper is organized as follows. In \cref{sec:prelim} we introduce our notation and present some preliminary results.  In \cref{sec:stability} we apply known stability   theorems for exponential systems to prove stability results for Fourier matrices. In \cref{sec:bounds} we use explicit estimates for singular values of  Fourier matrices to produce corresponding bounds of exponential systems. We summarize our conclusions and plans for future work in \cref{sec:future}. Proofs are contained in \cref{sec:proofframe,proofs}.

We believe that our stability results for Fourier matrices, derived through connections to exponential systems, will be valuable in numerical applications. The condition number of Fourier and Vandermonde matrices have been heavily studied in the literature and there are numerous applications. The condition number provides insight into the stability of solving a system of equations, like in trigonometric interpolation, where there are perturbations in the sampled values (i.e., in the range $y$). On the other hand, stability of nodes/frequencies deals with the sensitivity of Fourier matrices to sampling errors (i.e., in space $x$ or frequency $\xi$). 

There are few results for stability of Fourier matrices. A recent paper by A. Yu and A. Townsend \cite{YT} derived such a theorem for discrete Fourier transform matrices in the single variable case, equivalent to \cref{T3-freqstability} part (a) for $d=1$. We refer the reader to this article for further numerical applications and survey of this problem in the numerical analysis community. 
\subsection{Related results}

While completing this manuscript we became aware that T. Alemany and S. Nitzan were studying the lower Riesz basis bounds for exponential systems over an interval of the real line \cite{AN}. Although their main focus of interest  is considerably different than ours,  they independently obtained \cref{thm:kadec} and \cref{T-rectangular} for square matrices.

\medskip
In our paper we consider exponential bases and frames   on  finite union of hyper-cubes in $\R^n$   whose  vertices are in a lattice, and we exploit known results for these bases to prove corresponding results  for Fourier matrices.
Exponential bases and frames are known to exist on any finite union of hyper-cubes in $\R^n$ with sides parallel to the coordinate axes, 
and also on some infinite ones  (see e.g \cite{ AsDr,  KoNi, Marzo}).   In this general case, such exponential bases and their frame constants are often not explicit, which makes them less suitable for other applications.

\medskip

\section{Preliminaries}
\label{sec:prelim}

\subsection{General notation}\label{2.1}

We let $\N$ denote the natural numbers. For $\vec v=(v_1,\, ..,\, v_d)$,  and $\vec w= (w_1,\, ..,\, w_d)\in\R^d$,  we let $\vec v\cdot \vec w= \sum_{k=1}^d v_k w_k$ be their inner product and $\vec v \vec w=(v_1w_1,\dots,v_dw_d)$ be their component-wise multiplication. If the components of $\vec v$ are nonzero, we let ${\vec v\,}^{-1}= (v_1^{-1}, \dots, v_d^{-1})$ and $\vec w \,{\vec v\,}^{-1}=\vec w/\vec v = {\vec v\,}^{-1}\,\vec w = (\frac{w_1}{v_1},\, ... \frac{w_d}{v_d})$ be their component-wise quotient. For a given $\Omega\subset\R^d$ and  $\vec M\in\R^d$, we let  $\vec M \Omega:=\{\vec M \vec\omega\colon \vec \omega \in \Omega\}$ be  the dilation of $\Omega$ by $\vec M$. If   $\vec M\in\R^d$ with $M_k>0$ for all $k$, we define the rectangle $R(\vec M)= [0, M_1)\times \cdots \times [0, M_d)$.

As customary, the singular values of a matrix $\bA\in\C^{m\times n}$ are denoted by $0\leq \sigma_n(\bA)\leq \cdots \leq \sigma_1(\bA)$. When the dimensions of $\bA$ are unspecified, we will simply let $\sigma_{\max}$ and $\sigma_{\min}$ denote its biggest and smallest singular values (which could be zero). The ratio $\kappa(\bA):=\sigma_1/\sigma_n$ is called the {\it condition number} of $\bA$. 

For $p\ge 1$, we  use the notation $\|\cdot\|_p$   for the $L^p$ norm of a measurable function on a  domain of $\R^d$,  and  also  for the $\ell^p$ norm of a vector $x\in\C^m$, (that is,   $\|x\|_p=  (\sum_{j=1}^m |x_j|^p )^{\frac 1p})$ when there is no risk of confusion.  
 
%and the $L^p$ to $L^p$ (or $\ell^p$ to $\ell^p$) operator norm. 

\subsection{Exponential systems} We have used    the excellent textbook \cite{Cr}  for  the definitions  and properties of bases and frames. Let $H$ be a separable Hilbert space. We say a sequence of vectors $\mathcal{V}= \{v_j\}_{j\in\N}$ in $H$ is a {\it Riesz sequence} in $H$ if there are $A,B>0$ such that for all finite sequences $ \{a_j\}_{j\in J}\subset\C $, 
\begin{equation}\label{e2- Riesz-sequence}
	A \, \sum_{j\in J}   |a_j|^2   
	\leq  \Big\| \sum_{j\in J}  a_j  v_j \Big\|_H^2  
	\leq B \, \sum_{j\in J} |a_j|^2. 
\end{equation}
A sequence of vectors $\mathcal{V}= \{v_j\}_{j\in\N}$ in $H$ is a {\it frame} for $H$ if 
there exist  constants $A, \, B>0$    such that for  every $w\in H$,
\begin{equation}\label{e2-frame}
	A \, \|w\|^2_H
	\leq  \sum_{j=1}^\infty |\l  w, v_j\r_H |^2
	\leq B \, \|w\|_H^2.
\end{equation}
Here $\langle \cdot, \cdot \rangle_H$  and $\|\cdot \|_H=\sqrt{\l \cdot , \cdot \r_H} $  are the  inner product  and norm  in $H$. If $\calV$ is both a Riesz sequence and frame, then we say it is a {\it Riesz basis}.

Generally, $A$ and $B$ should be interpreted as any valid bounds for which \eqref{e2- Riesz-sequence} holds when dealing with Riesz sequences, or when \eqref{e2-frame} is satisfied for frames. The optimal constants are the largest $A$ and smallest $B$ for  which the inequalities above hold. With some abuse of terminology, we term  $\kappa:=\sqrt{B/A}$      the  {\it condition  number}  of the system ${\mathcal V}$.   We will mainly consider $H=L^2(D)$, with $D\subset \R^d$ of finite Lebesgue measure $|D|$. For a discrete $\Lambda\subset \R^d $, we call $\B(\Lambda)=\{ e^{2\pi i \vec\lambda \cdot x}\}_{ \vec\lambda\in\Lambda}$ an {\it exponential system}. %on $D$.  

\subsection{Fourier matrices} \label{2.3} In this paper, a (generalized) {\it Fourier matrix} refers to any type of matrix $\bF$ for which there are finite sets  $\Omega\subset\Z^d$ and $U\subset (\R/\Z)^d={\mathbb T}^d$ such that
\begin{equation}
	\label{e-Fouriermatrix}
	\bF:= \bF(\Omega,U):= \Big[ \, e^{2\pi i\vec \omega\cdot \vec u} \, \Big]_{\vec \omega\in \Omega\atop{ \vec u\in U}}. 
\end{equation}
We refer to $u\in U$ as a {\it node} and $U$ as the {\it node set}, whereas $\Omega$ is the {\it frequency set}.

Our definition of a Fourier matrix includes numerous classes of matrices that have been carefully studied. Recall that a {\it Vandermonde matrix} $\bV$ of size $L\times N$  with nodes on the complex unit circle $\T=\R/\Z$ is of the form
$$
\bV:=\bV(L,X):=\Big[ \, e^{2\pi i jx_k}\,  \Big]_{0\leq j<L\atop{ 1\leq k\leq N}}, \wherespace X=\{x_k\}_{k=1}^N \subset \T. 
$$	
It is well known that $\bV$ has rank $N$ if and only if  $X$ consists of $N$ distinct elements of $\T$. 

Our definition of a Fourier matrix also includes the usual (unnormalized) {\it discrete Fourier transform} (DFT) matrix. We recall that $\bF(\Omega,U)$ is a DFT matrix if for some $\vec M\in \N^d$,
$$
\Omega=R(\vec M)\cap\Z^d \andspace U=\Big\{\frac {\vec k}{\vec M} \colon 0\leq k_j <M_j\text{ for each } j \Big\}.
$$
When $U'$ is a proper subset of $U$, then $\bF(\Omega,U')$ is  a sub-matrix of the DFT matrix  $\bF(\Omega,U)$ formed by selecting a subset of the latter's columns. Note that partial DFT matrices have orthogonal columns while general Fourier matrices do not.  

\subsection{Exponential systems and Fourier matrices} 

Fourier matrices naturally appear in our investigation of exponential systems  via \cref{T1-frame}. The most basic connection is given by the following example. For a given $\vec M\in\N^d$, consider the exponential basis $\B:=\{e^{2\pi i \vec M^{-1} \vec n\cdot x}\}_{\vec n\in \Z^d}$ on  the rectangle $R(\vec M)$  defined in the  sub-section \ref{2.1}. Notice that 
$$
\B=\bigcup_{j_1=0}^{M_1-1} \cdots \bigcup_{j_d=0}^{M_d-1}\big\{ e^{2\pi i (\vec n+\frac {\vec j}{\vec M}) \cdot x}\big\}_{\vec n\in\Z^d}. 
$$

Letting $P=R(\vec M)\cap\Z^d$ and $\Delta=\{\frac {\vec j}{\vec M}\colon \vec j\in P\}$, notice that $R(\vec M)=T(P)$ and $\B=\B(\Delta)$. Hence, the matrix $\G(\Delta,P)$ associated with the standard exponential basis on a rectangle with integer side lengths is just a DFT matrix. For example, when   $d=1$,  we have 
\begin{equation}\label{e-st-FM}
	\G(\Delta, P)= \Big[e^{2\pi i jk/M}\Big]_{0\leq k\leq N-1\atop{0\leq j\leq L}}.
\end{equation}
To relate $\G$ and $\bF$, we have the most immediate connection which is that $\G(\Delta,P)=\bF(P,\Delta)$. However, there are circumstances where it is advantageous to construct $\Omega$ and $U$ that are different from $P$ and $\Delta$ such that $\G(\Delta,P)=\bF(\Omega,U)$. 

For instance, whenever $\Delta$ or $P$ is an arithmetic progression, $\G(\Delta,P)$ is a Vandermonde matrix. Indeed, if $\vec \delta_j=(j-1)\vec\delta$ for $j=1,\dots,L$ and some $\vec \delta\in\R^d$, then $\G(\Delta,P)=\bV(L,X),$ where $X=\{\vec \delta \cdot \vec p_k\}_{k=1}^N$ are the nodes. On the other hand, if $\vec p_k = (k-1)\vec p$ for $k=1,\dots,N$ and $\vec p\in\Z^d \setminus\{0\}$, then $\G(\Delta, P)^*=\bV(N,Z),$ where $X=\{\vec \delta_j\cdot \vec p\}_{j=1}^L$ are the nodes. 

Vandermonde matrices are inherently one-dimensional. To obtain a multivariate Fourier matrix from $\G(\Delta,P)$, we can consider the case  when $\Delta$ are scaled integers. Fix some $\vec M\in \N^d$ and %recall the rectangle $R(\vec M)$ from our earlier notation. 
set $\Omega:=R(\vec M)\cap\Z^d$ and $\Delta:=\vec M^{-1}\Omega$. For any $P\subset \Omega$, we let $U:=\vec M^{-1} P$. Then we have $\G(\Delta,P)=\bF(\Omega,U).$ In particular, this is a DFT matrix if $P=\Omega$ and is a partial DFT matrix if $P$ is a proper subset of $\Omega$. More generally, suppose $\Omega\subset\Z^d$ is an arbitrary set (e.g., common examples are multi-integers in a cube or ball). For any $\alpha >0$ and $P\subset \Z^d$, let $\Delta=\alpha \Omega$ and $U=\alpha P$. Then we have $\G(\Delta,P) = \bF(\Omega,U).$

By using the above observations together with \cref{T1-frame}, the following result generalizes Theorem 1.3 in \cite{DC}. Its proof is contained in \cref{proofs}.

\medskip\begin{Thm}\label{T1-special-delta}
	Let $\vec\delta\in \R^d$ be fixed and $P=\{\vec p_k\}_{k=1}^N\subset\Z^d$. For any $L\geq N$, the set  
	\begin{equation}\label{d2-def-S}
		\B
		=\bigcup_{j=1}^{L} \big\{e^{ 2\pi i  (\vec n +(j-1)\vec \delta) \cdot x}\big\}_{\vec n\in\Z^d}   
	\end{equation} 
	is a frame for $L^2(T(P))$ if and only if 
	\begin{equation}\label{e2-cond=delta} \l \vec p_k - \vec p_{k'},\ \vec \delta\r\not\in\Z,
		\quad   \mbox{ for every $k, k' \in [1,\,  N]$,   $k\neq k'$}.
	\end{equation} 	
	The optimal frame constants   of $\B$  are  the squares of the maximum and minimum singular value of the  $L\times N$  Vandermonde matrix $\bV(L,\Omega) $. 
\end{Thm}

\subsection{Kadec theorem and generalizations}
 
Riesz bases are {stable}, in the sense that any sufficiently small perturbation of a Riesz basis produces a Riesz basis. The problem is to quantify  what is  the largest allowable perturbation, or the {\it stability bound} of the basis. 

A workhorse result in this direction is the Paley-Wiener theorem and its variations,   which can  be used to prove stability results for bases on  general Hilbert spaces. See \cite{PW} and \cite{O}. A celebrated stability result for  the standard orthogonal basis $\{e^{2\pi i n x}\}_{n\in\Z}$ on $L^2(0,1)$ is the classical Kadec theorem; see \cite{K} and also \cite[Chaper 1, Theorem 14, page 42]{Y}. Kadec's theorem has been generalized, and we will state a few results momentarily, after introducing some additional notation.  

We say a nontrivial sequence $\vec a\colon \Z^d\to\R^d$ is (at most) {\it rank one} if there exist $a_1,\dots,a_d\colon \Z\to\R$ such that for each $\vec n\in \Z^d$,
$$
\vec a(\vec n)=(a_1(n_1),\dots,a_d(n_d)).
$$ 
The main implication of this definition is that only one component of $\vec a$ changes when the corresponding entry of $\vec n$ changes. The rank one terminology will make sense once we deal with matrices.  We define the functions 
\begin{align}
	C(t) &:=1-\cos(\pi t)+\sin(\pi t)=1-\sqrt 2 \sin\left(\frac \pi 4 - \pi t\right), \label{eq:Cconst} \\
	D(t) &:= \left(1-\cos(\pi t)+\sin(  \pi t) +\frac{\sin(\pi t)}{\pi t}\right)^d-\left(\frac{\sin(\pi t)}{\pi t}\right)^d. \label{eq:Dconst}
\end{align}

\begin{Thm} 
	\label{thm:kadec} 
	Let $\vec \lambda\colon \Z^d\to\R^d$ and suppose $ \{e^{2\pi i  \vec \lambda({\vec n}) \cdot  x}\}_{{\vec n}\in\Z^d}$ is a Riesz basis for (or a Riesz sequence in, or a frame for) $L^2(Q)$ with bounds $A$ and $B$. Let $\vec \delta\colon \Z^d\to\R^d$ such that $\ell:=\|\vec\delta\|_\infty<\frac 14.$	
	\begin{enumerate}[{\rm (a)}]		
		\item
		If $D(\ell)<\sqrt{\frac AB},$ then 
		$\{e^{2\pi i (\vec\lambda(\vec n)+\vec\delta(\vec n))\cdot x}\}_{\vec n\in\Z^d}$ is a Riesz basis for (or a Riesz sequence in, or a frame for) $L^2(Q)$  with bounds
		$$ 
		A'=	A \left(1-\sqrt{\frac BA} \, D(\ell) \right)^{2} \andspace B'= B (1+ D(\ell) )^2.
		$$ 
		\item 
		If $\lambda$ and $\delta$ are rank one, and $C(\ell)<\sqrt{\frac AB},$ then $\{e^{2\pi i (\vec\lambda(\vec n)+\vec\delta(\vec n))\cdot x}\}_{\vec n\in\Z^d}$ is a Riesz basis for (or a Riesz sequence in, or a frame for) $L^2(Q)$ with bounds
		\begin{equation*}
			A'=	A \left(1-\sqrt{\frac BA} \, C( \ell)\right)^{2d} \andspace 
			B'= B (1+ C( \ell))^{2d}.
		\end{equation*}
	\end{enumerate}
\end{Thm} 

\cref{thm:kadec} part (a) is the general case, which does not require that $\vec \lambda$ and $\vec \delta$ are rank one, and was proved  by W. Sun  and X. Zhou (see \cite[Theorem 1.3]{SZ}). Part (a) for $d=1$ was proved by R. Balan \cite{B}  in  \cref{thm:kadec}. Part (b) is for the special case where both $\vec \lambda$ and $\vec \delta$ are rank one and was shown in \cite[Thm 1.2]{SZ}. Technically,  the results in  \cite{B,SZ} are for frames, but an inspection of the  proofs shows that the same strategies carry over to the Riesz sequence and basis cases. We also remark that $A'$ and $B'$ are valid, not necessarily optimal, constants. Also, when $d=1$, the rank one assumption is trivially satisfied and $D(\ell)=C(\ell)$, so parts (a) and (b) agree.

\begin{Rem} \label{Rem-kadec}
	Let $\vec M\in \N^d$ and $R:=R(\vec M)$. It is easy to verify, using a change of variables and prior notation that $ \{e^{2\pi i \vec M^{-1} (\vec n+\vec\delta(\vec n))\cdot  x}\}_{{\vec n}\in\Z^d}$ is a Riesz basis for $L^2(R)$ whenever the assumptions of \cref{thm:kadec} are satisfied by $\ell=\|\vec\delta\|_\infty$. %In which case
	The constants of the  basis 
	 $ \{e^{2\pi i \vec n/\vec M    \cdot  x}\}_{{\vec n}\in\Z^d}$ are  $A=B=M_1\cdots M_d$ and valid constants $A'$ and $B'$ can be explicitly computed in terms of the quantities $D(\ell)$ and $C(\ell)$   given in \eqref{eq:Cconst} and \eqref{eq:Dconst}, respectively.
\end{Rem}

In the classical Kadec theorem, it is known that $\frac 1 4$ cannot be replaced with a bigger number. By considering more structured perturbations, the $\frac 1 4$ barrier can be overcome. The following result  is a  higher-dimensional version of  \cite[Corollary 5.1] {DC}. Its proof is in Section \ref{proofs}.
 
\begin{Thm}\label{T-New-Kadec}
	Let $\vec M\in \N^d$, $J=R(\vec M)\cap \Z^d$, and $\vec \epsilon:J\to \R^d$ be rank one. The set
	\begin{equation*} 
		\B
		= \bigcup_{\vec j\in J} \big\{e^{ 2\pi i(\vec n + (\vec j+\vec \epsilon(\vec j))/\vec M) \cdot \vec x} \big\}_{\vec n \in\Z^d}  
	\end{equation*}
	is a Riesz basis for $L^2(R(\vec M))$ if and only if, for every $\vec i,\, \vec j\in J$, with $\vec i\ne \vec j$  and every $k=1,\dots, d$, we have that  
	\begin{equation} \label{e-cond1}
		\frac 1 {M_k} \, \big(j_k-  i_k+\epsilon_k(j_k)-\epsilon_k(i_k) \big) \not\in\Z.
	\end{equation}  
	In particular, ${\B}$ is a Riesz basis for $L^2(R(\vec M))$  whenever $\|\vec\epsilon\|_\infty< {\frac 12}$.
\end{Thm}

  The frame constants of the basis ${\B}$ cannot easily be evaluated when $\frac 14\leq \|\vec \epsilon\|_\infty< \frac 12$. Under an assumption on averaged perturbations, Avdonin \cite{A} proved a powerful stability theorem,   but the frame constants of the perturbed basis are not explicit. Explicit  frame bounds in   Avdonin's theorem  can be found in    \cite{AN}.

\section{Stability of Fourier matrices}

\label{sec:stability}

We first examine what happens when the frequencies $\Omega$ of a given Fourier matrix are perturbed while the nodes $U$ are fixed. Here, $\Omega\subset\Z^d$, and its perturbation is denoted  $\Omega'\subset \R^d$. By slight abuse of notation, we also let $\bF(\Omega',U)$ be the matrix in \eqref{e-Fouriermatrix} with $\Omega'\subset\R^d$ not necessarily a subset of the integers and $U\subset [0, 1]^d$. 

We begin by deriving a basic but loose estimate via tools from linear algebra and matrix perturbation theory, which will serve as a benchmark comparison. Its proof is in \cref{proofs}. 

\begin{Prop} \label{prop:badstability2}
	Let $\Omega=\{\vec \omega_k\}_{k=1}^L\subset \Z^d$ and $U\subset\T^d$,    $|U|=N$. For any $p\in [1,\infty]$ and $\epsilon>0$, let $\Omega':=\{\vec \omega_k'\}_{k=1}^L$ such that $\|\vec \omega_k-\vec \omega_k'\|_p \leq \epsilon$ for each $k\in \{1,\dots,L\}$. Then 
	\begin{align*}
		\sigma_N(\bF(\Omega',U))
		&\geq \sigma_N(\bF(\Omega,U))-\pi d^{1/p'} \sqrt{L N} \epsilon, \\ 
		\sigma_1(\bF(\Omega',U))
		&\leq \sigma_1(\bF(\Omega,U)) + \pi d^{1/p'} \sqrt{L N} \epsilon. 
	\end{align*}
\end{Prop}  

We start with a common situation where the frequencies of a discrete Fourier transform matrix are perturbed. 

\begin{Thm}\label{T3-freqstability}
	For any $\vec M\in \N^d$, let $\Omega:=R(\vec M)\cap\Z^d$ and $U\subset \vec M^{-1} \Omega$. For any $\vec \epsilon:\Omega \to \R^d$, define $\Omega':= \{\vec j + \vec \epsilon(\vec j)\colon \vec j \in \Omega\}$ and $\ell:=\|\vec \epsilon\|_\infty$. 
	\begin{enumerate}[{\rm (a)}]
		\item 
		If $\ell < {\frac 14}$, then the extreme singular values of $\bF(\Omega',U)$ satisfy
		$$	
		\sigma_{\min}
		\ge (1- D(\ell) ) \, \sqrt{M_1\cdots M_d} \andspace \sigma_{\max}
		\leq (1+D(\ell) ) \, \sqrt{M_1\cdots M_d}. 
		$$ 
		\item 
		If $\ell < {\frac 14}$ and $\vec \epsilon$ has rank one, then the extreme singular values of $\bF(\Omega', U)$ satisfy
		$$
		\sigma_{\min}
		\ge (1- C(\ell))^d \, \sqrt{M_1\cdots M_d} \andspace \sigma_{\max}
		\leq (1+C(\ell))^d \, \sqrt{M_1\cdots M_d}. 
		$$ 
		\item 
		If $\vec \epsilon$ has rank one, then
		$\bF(\Omega', U)$ is nonsingular if and only if 
		for every $\vec i, \, \vec j\in \Omega$, with $i_k \ne j_k$ for every $k=1,\dots,d$,  we have that  \begin{equation*}
			\frac 1 {M_k} \big(j_k-i_k+\epsilon_k(j_k)-\epsilon_k( i_k) \big)  \not\in\Z.
		\end{equation*}  
	\end{enumerate}	 
\end{Thm}

\begin{proof}
	Let $U_0=\vec M^{-1} \Omega$. Since $U\subset U_0$, we see that $\bF(\Omega',U)$ is a submatrix of $\bF(\Omega',U_0)$ formed by deleting columns from the latter. By the variational characterization of the smallest and largest singular values of a matrix, we see that 
	$$
	\sigma_{\min}(\bF(\Omega',U)) \geq \sigma_{\min}(\bF(\Omega',U_0)) 
	\andspace \sigma_{\max}(\bF(\Omega',U)) \leq \sigma_{\max}(\bF(\Omega',U_0)). 
	$$
	Thus, from now on, we assume that $U=U_0$.   
	
	First, we define $\Delta:=\vec M^{-1}\Omega$ and $P=\Omega$ so that $T(P)=R(\vec M)$. Notice that $\bF(\Omega,U)=\G(\Delta,P)$ and the associated exponential system is 
	$$
	\B(\Delta)
	=\bigcup_{\vec j\in \Omega} \big\{e^{2\pi i (\vec n+\frac {\vec j}{\vec M}) \cdot x} \big\}_{\vec n \in\Z^d} 
	=\big\{e^{2\pi i \frac {\vec n}{\vec M}\cdot x} \big\}_{\vec n \in\Z^d}.
	$$
	This is the standard orthogonal basis for $L^2(R(\vec M))$ with optimal lower and upper constants $A=B=M_1\dots M_d$. 
	
	Next, we define $\Delta'=\vec M^{-1}\Omega'$ and extend $\vec \epsilon$ to a $\Omega$ periodic function on $\Z^d$. \cref{thm:kadec}and the change of variables argument outlined in \cref{Rem-kadec} provide sufficient conditions on $\vec \epsilon$ such that the perturbed family 
	$
	\{e^{2\pi i \vec M^{-1}(\vec n+\vec \epsilon)\cdot x} \}_{\vec n \in\Z^d}
	$
	is an exponential basis or frame for $L^2(R(\vec M))$. %To put this into our framework, 
	Indeed, observe that
	$$
	\B(\Delta')	
	=\bigcup_{\vec j\in \Omega} \big\{e^{2\pi i (\vec n+\frac{\vec j}{\vec M} + \frac {\vec \epsilon}{\vec M})\cdot x} \big\}_{\vec n \in\Z^d}
	=\big\{e^{2\pi i \vec M^{-1}(\vec n+\vec \epsilon)\cdot x} \big\}_{\vec n \in\Z^d}. 
	$$
	By \cref{T1-frame} part (a), $\B(\Delta')$ is an exponential frame for $L^2(R(\vec M))$ with constants $A'$ and $B'$ if and only if $\G(\Delta',P)$ has extreme singular values $\sqrt{A'}$ and $\sqrt{B'}$. Finally, a calculation shows that $\G(\Delta',P)=\bF(\Omega',U)$. Parts (a) and (b) follow from the corresponding parts of \cref{thm:kadec}, while part (c) follows from \cref{T-New-Kadec}.   
\end{proof}

\begin{Rem}
	Let us compare \cref{prop:badstability2} with  \cref{T3-freqstability} when
	$\Omega=R(\vec M)\cap \Z^d$ and $U=\vec M^{-1}\Omega.$ Then $\bF(\Omega,U)$ is the DFT matrix and has extreme singular values both equal to $\sqrt{M_1\cdots M_d}$. We ask, what is the largest allowable perturbation of $\Omega$ to $\Omega'$ in $\ell^\infty$ for which  %such that the resulting $\bF(\Omega',U)$ enjoys 
	\begin{equation}
		\label{eq:desired}
		\sigma_{\min}(\bF(\Omega',U))
		\geq \frac 1 2 \, \sigma_{\min}(\bF(\Omega,U))? 
	\end{equation}
	While the $\frac 1 2$ factor was chosen arbitrarily, it is natural to investigate a ratio of the perturbed versus unperturbed singular values. 
	
	We use \cref{prop:badstability2} with $p=\infty$, $L=N=M_1\cdots M_d$, and $p'=1$. To get the desired inequality \eqref{eq:desired} from this proposition, we need to set 
	$$
	\epsilon 
	= \frac{1}{2\pi d \sqrt{M_1\cdots M_d}}. 
	$$
	
	We use \cref{T3-freqstability}. Call $\ell_a$ and $\ell_b$ the largest $\ell\in [0,\frac 1 4)$ for which $1-D(\ell) = \frac 12$ and $(1-C(\ell))^d = \frac 1 2$, respectively. To use this theorem, the perturbation $\vec \epsilon$ must satisfy 
	$\|\vec \epsilon\|_\infty<\ell_a$ and $\|\vec \epsilon\|_\infty<\ell_b$ respectively, and part (b) further requires $\vec\epsilon$ to be rank one. 
	
	For fixed dimension, \cref{T3-freqstability} is clearly superior to \cref{prop:badstability2} because the set of allowable perturbations does not shrink as the entries of $\vec M$ increase. Let us briefly comment on the fairness of this comparison. \cref{prop:badstability2} is more general, since it holds for any $\Omega$ and $U$, while \cref{T3-freqstability} only holds for the DFT matrix. On the other hand, this is also a favorable setting for \cref{prop:badstability2}, since the allowable perturbation size is dependent on the extreme singular values of the original Fourier matrix, which are equal when it has orthogonal columns.
\end{Rem}

\cref{T3-freqstability} part (a) for $d=1$ was established in a recent publication of A. Yu and A. Townsend \cite[Corollary 2.3]{YT}. Their proof closely follows a standard proof of Kadec's theorem, while our proof relies on a general connection between exponential systems and Fourier matrices given in \cref{T1-frame}. Through this connection, we obtained a different proof. To our best knowledge, \cref{T3-freqstability} part (a) for $d>1$, part (b) and part (c) are completely new. 

At this point, it is hopefully clear why we used the rank one definition from before. In \cref{T3-freqstability} parts (b) and (c), the assumption that $\vec \epsilon$ has rank one means  that the matrix $\bF(\Omega,U)$ is perturbed to $\bF(\Omega',U)$ by a matrix $\bE$ that depends on $\vec\epsilon$. It has rank one since $\vec\epsilon$ has rank one. 

Next, we examine what happens when the nodes $U$ are perturbed   while the frequencies $\Omega$ are fixed. Let $U =\{\vec u_k\}_{k=1}^N$ and let  $U'=\{\vec u_k'\}_{k=1}^N$ be the perturbed nodes. As usual, we let $\Omega=\{\vec \omega_j\}_{j=1}^L$. We start by deriving a loose estimate of the singular values of $\bF(\Omega,U')$ based on general tools that will serve as a benchmark comparison.

\begin{Prop} \label{prop:badstability}
	For any $N,L\in \N$, $\Omega\subset \Z^d$, $p\in [1,\infty]$, and $\epsilon>0$ such that $\|\vec u_k- \vec u_k'\|_p\leq \epsilon$ for each $k\in \{1,\dots,N\}$, letting $C_{\Omega,p}:= 2\pi \sqrt{\sum_{j=1}^L \|\vec \omega_j\|_{p'}^2}$, we have  
	\begin{align*}
		\sigma_N(\bF(\Omega,U'))
		&\geq \sigma_N(\bF(\Omega,U))-C_{\Omega,p} \sqrt N \epsilon, \\ 
		\sigma_1(\bF(\Omega,U'))
		&\leq \sigma_1(\bF(\Omega,U)) + C_{\Omega,p} \sqrt N \epsilon. 
	\end{align*}
\end{Prop}

This proposition yields a nontrivial lower bound for $\sigma_N(\bF(\Omega,U'))$ only when $\epsilon$ is sufficiently small. This is a rather disappointing result since $\epsilon$ is required to become progressively smaller as the matrices' size increase. A primary reason for this deficiency is that this result hardly uses any properties of Fourier matrices. 

It is natural to derive a result for node stability using Kadec theorems for exponential frames. The following is a one-dimensional version for Vandermonde matrices.

\begin{Thm} \label{T-rectangular}
	Let $U=\{u_k\}_{k=1}^N\subset \R$ and assume the Vandermonde matrix $\bV:=\bV(L,U)$ has rank $r=\min\{L, N\}$. Let $\{\delta_k\}_{k=1}^N \subset \R$ such that $\ell:=\|\delta\|_\infty$ satisfies $C(\ell)<\frac{\sigma_r(\bV)}{\sigma_1(\bV)}$. 
	Letting $U'=\{u_k+\frac{\delta_k}{L}\}_{k=1}^N$, the $L\times N$ matrix $\bV':=\bV(L,U')$ has rank $r$, and its extreme singular values satisfy the inequalities
	\begin{equation}\label{e-new-ineq}
		\sigma_r(\bV') \ge \sigma_r(\bV) \left(1- \frac{\sigma_1(\bV)}{\sigma_r(\bV)} \,C(\ell)\right)
		\andspace
		\sigma_1(\bV') \leq \sigma_1(\bV) (1+C(\ell)).
	\end{equation}    
\end{Thm}

\begin{proof}
	We only prove this theorem for the case $N\leq L$, since the other case can be proved in a similar fashion. Since $\bV$ is a $L\times N$ matrix with rank $r=\min\{L,N\}=N$, we will use \cref{T1-frame} part (b) where $U$ corresponds to $\Delta$ and $\{0,\dots,L-1\}$ plays the the role of  $P$. By the theorem,
	$
	\bigcup_{k=1}^N \big\{e^{2 \pi i (n+u_k) x}\big\}_{n \in \Z} 
	$
	is a Riesz sequence in $L^2([0,L))$ with optimal bounds $A=\sigma_r^2(\bV)$ and $B=\sigma_1^2(\bV)$. 
	
	Namely, for all $f\in L^2([0,L))$, we have
	\begin{equation*}
		A \|f\|_{L^2([0,L))}^2 
		\leq \sum_{k=1}^N \sum_{n\in \Z} \left| \int_0^L e^{2\pi i (n+u_k)x} f(x)\, dx\right|^2
		\leq B \|f\|_{L^2([0,L))}^2.
	\end{equation*}
	We let $f_L(x):=f(Lx)$ so that $f_L\in L^2([0,1))$. After a change of variable in the $L^2$ norm, we obtain $\|f\|_{L^2([0,L))}^2=L \|f_L\|_{L^2([0,1))}^2$, and 
	$$
	\left| \int_0^L e^{2\pi i (n+u_k)x} f(x)\, dx\right|^2
	= \left| \int_0^1 e^{2\pi i L (n+u_k)t} f(Lt) L \, dt\right|^2
	= L^2 \left| \int_0^1 e^{2\pi i L (n+u_k)t} f_L(t) \, dt\right|^2. 
	$$
	Combining the above, we obtain
	\begin{equation*}
		\frac AL \|f_L\|_{L^2([0,1))}^2 
		\leq \sum_{k=1}^N \sum_{n\in \Z}
		\left| \int_0^1 e^{2\pi i L (n+u_k)t} f_L(t) \, dt\right|^2
		\leq \frac B L \|f_L\|_{L^2([0,1))}^2.
	\end{equation*}
	Note that $f\mapsto f_L$ is a bijection between $L^2([0,L))$ and $L^2([0,1))$. This now shows that 
	$$
	\bigcup_{k=1}^N \big\{e^{2 \pi i L(n+u_k) x}\big\}_{n \in \Z},
	$$
    is a Riesz sequence in $L^2([0,1))$ with optimal bounds $L^{-1} A$ and $L^{-1}B$. By \cref{thm:kadec} part (a), whose use is justified by the assumption that $C(\ell)< \sqrt{A/B}$, we deduce that 
    $$
    \bigcup_{k=1}^N \big\{e^{2 \pi i (Ln+Lu_k+\delta_k) x}\big\}_{n \in \Z},
    $$ is a Riesz sequence in $L^2([0,1))$ with bounds $L^{-1} A'$ and $L^{-1} B'$, where $A':=A\big(1-\sqrt{\frac BA} C(\ell)\big)^2$ and $B':=B(1+C(\ell))^2$. After  undoing this scaling we can conclude that $\bigcup_{k=1}^N \big\{e^{2 \pi i (n+u_k+\delta_k/L) x}\big\}_{n \in \Z}$ is a Riesz sequence in $L^2([0,L))$ with bounds $A'$ and $B'$. To complete the proof, we use \cref{T1-frame} part (b) again to see that the corresponding matrix of this Riesz sequence is $\bV(L,U')$.
\end{proof}

Before we move on, let us explain the meaning of this theorem. If the original Vandermonde matrix $\bV$ has good condition number, then as expected, we are allowed bigger perturbations, as seen in the requirement that $C(\ell)<{\sigma_r(\bV)}/{\sigma_1(\bV)}=\kappa^{-1}(\bV)$. The perturbed matrix has nodes $U'=\{u_k+\frac{\delta_k}L\}_{k=1}^N$, so while $\|\delta\|_\infty<\frac 1 4$, the nodes are perturbed by at most $\frac 1 {4L}$ with respect to the wrap-around metric  $|t-s|_\T:=\min_{n\in \Z}|t-s-n|.$ This $\frac 1 L$ factor appears because $L$ corresponds to the number of frequencies in $\bV$ so we expect perturbations of the nodes to be inversely proportional to $L$. 

\begin{Rem}
	To compare \cref{prop:badstability} and \cref{T-rectangular}, we seek to find the set of allowable perturbations such that $\sigma_r(\bV')\geq \frac 1 2 \sigma_r(\bV).$ If we use \cref{prop:badstability} with $\Omega=\{0, 1,... L-1\}$  and  $p=\infty$, the maximum allowable perturbations have size 
	$$
	\epsilon
	\leq  \frac{\sigma_r(\bV)}{2C_{\Omega,1} \sqrt N}
	= \frac{\sigma_r(\bV)}{4 \pi \sqrt{N\sum_{k=0}^{L-1} k^2}}.
	$$	
	The denominator is $O(N^{1/2} L^{3/2})$. On the other hand, \cref{T-rectangular} allows the nodes to be perturbed by $\frac \ell L$ provided that 
	$ 
	C(\ell)< \frac{\sigma_r(\bV)}{2\sigma_1(\bV)}. 
	$ 
	To obtain an explicit comparison, notice that $\sigma_1(\bV)\in [\sqrt L, \sqrt{LN}]$ since each column of $\bV$ has $\ell^2$ norm $\sqrt L$ and there are $N$ columns. Thus, this shows that \cref{T-rectangular} always requires a weaker assumption than \cref{prop:badstability} whenever $L$ is large  enough. 
\end{Rem}

A version of \cref{T-rectangular} for multivariate Fourier matrices   can derived using analogous arguments, but we  do not include it here. 

\section{Explicit bounds to exponential systems}

\label{sec:bounds}
 
Fourier matrices have been studied quite extensively by numerical and harmonic analysts due to their close relation to polynomial interpolation and approximation \cite{Z}, the mathematics of super-resolution \cite{liliao}, compressed sensing \cite{FR}, and more. Accurate estimates for the condition number of such matrices are useful for determining the stability of numerical algorithms to noise and other perturbations. 

The condition number of Vandermonde matrices greatly depends on  the pairwise distance of the nodes  according to the wrap-around metric  $|t-s|_\T:=\min_{n\in \Z}|t-s-n|.$ For a finite set $U=\{u_k\}_{k=1}^N\subset\R$, we define the minimum separation   
$$
\sep(U):= \min_{j\not=k} |u_j-u_k|_\T. 
$$
The following theorem was established in \cite{AB}. 

\begin{Thm}
	\label{thm:wellsep}
	For any finite set $U\subset\R$ and integer $L\geq 1$, if $\sep(U)>\frac 1 L$, then 
	$$
	L-\frac 1{\sep(U)}
	\leq \sigma_{\min}^2(\bV) \leq \sigma_{\max}^2(\bV)\leq L +\frac 1{\sep(U)}. 
	$$
\end{Thm}

Note that this theorem only applies to the setting where $L> |U|$, otherwise the condition $\sep(U)>\frac 1 L$ cannot be fulfilled. Due to \cref{T1-frame}, there is an analogous result for Riesz sequences. 

\begin{Cor}
 	Let $L\in\N$ and $\Delta=\{\delta_j\}_{j=1}^N\subset\R$ such that $\sep(\Delta)>\frac 1L$. Then $\bigcup_{k=1}^N \{e^{2\pi i (n+\delta_k)x}\}_{n\in \Z}$ is a Riesz sequence in $L^2([0,L))$ with bounds satisfying 
	$$
	A \geq L-\frac 1{\sep(\Delta)} \andspace B\leq L +\frac 1{\sep(\Delta)}.
	$$
\end{Cor}

Notice that \cref{thm:wellsep} becomes vacuous when  $\sep(U)\to \frac 1 L$. This is necessarily an artificial issue because $\bV$ is a Vandermonde matrix and will always have positive smallest singular value whenever $L\geq |U|$. For the case where the minimum separation is arbitrarily small, we recall \cite[Corollary 1 and Theorem 3]{limultiscale}. In the following, ${\rm diam}(U)$ is the diameter of $U$ and $\dist(U,V)$ is the distance between $U$ and $V$ with respect to $|\cdot|_\T$. 

\begin{Thm}\label{thm:li23}
	Let $L,N\in\N$ such that $L\geq 6N$. For arbitrary $\alpha \in (0,\frac 1 L)$ and $\lambda\in \N$, let $U\subset\R$ be a set of cardinality $N$ such that $\sep(U)\geq \alpha$ and there is a $r\in \N$ and disjoint union $U=U_1\cup \cdots \cup U_r$ such that $\max_k |U_k|\leq \lambda$. If $r\geq 2$, then also assume that $\beta:=\min_{j\not=k} \dist(U_j,U_k)\geq \frac {3\lambda}L$ and $\max_k {\rm diam}(U_k)< \beta$. Then there exist explicit universal constants $C,c>0$ such that 
	$$
	\sigma_N(\bV(L,U)) 
	\geq C \sqrt{\frac LN} \, (c L \alpha)^{\lambda-1} \andspace \sigma_1(\bV(L,U)) \leq \sqrt{L \big(\lambda+\frac 1 3 \big)}. 
	$$
\end{Thm}

A few comments are in order. The theorem roughly assumes that $U$ can be decomposed into subsets $U_1,\dots,U_r$, where each $U_k$ has cardinality at most $\lambda$, points within each $U_k$ cannot be too far apart from each other, and the distances between different $U_j$ and $U_k$ must be big enough.  The conclusion is a localization result -- when $\alpha$ is small, the dominant behavior of $\sigma_N(\bV(L,U))$ is $(L\alpha)^{\lambda-1}$, which only depends on how many elements of $U$ are contained within an interval of length $O(\frac 1 L)$ and their normalized separation $L\alpha$. This lower bound attains the optimal dependence on $\lambda,\alpha,L$, so the only potential looseness comes from the universal constants and behavior in $N$. Formulas for $C$ and $c$ as well as multiscale generalizations are available in \cite{limultiscale}. 

This theorem has a corresponding result for Riesz sequences. 

\begin{Cor}
	Let $L,N$ such that $L\geq 6N$. For arbitrary $\alpha \in (0,\frac 1 L)$ and $\lambda\in \N$, suppose $\Delta=\{\delta_k\}_{k=1}^N\subset\R$ satisfies the conditions of \cref{thm:li23} with $\Delta$ in place of $U$. There are universal constants $C,c>0$ such that $\bigcup_{k=1}^N \{e^{2\pi i (n+\delta_k)x}\}_{n\in \Z}$ is a Riesz sequence in $L^2([0,L))$ with bounds satisfying 
	$$
	A
	\geq \frac {CL}N \big(c L \alpha \big)^{2\lambda-2} \andspace 
	B \leq L \big(\lambda+\frac 1 3 \big). 
	$$
\end{Cor}

%Here, we have presented only a few representative results. There are estimates for the extreme singular values of random DFT matrices and higher dimensional examples. They have corresponding statements for Riesz sequences,  
There are estimates for the extreme singular values of random DFT matrices and higher dimensional examples. They have corresponding statements for Riesz sequences,   but here we  have presented only a few representative results. 

\section{Remarks and future work}
\label{sec:future}
 
Prior works have exploited connections between exponential systems and certain types of Fourier matrices. The purpose of our article is to provide a systematic study. We determined how choices of $\Delta$ and $P$ influences the singular values of  Fourier matrices $\bF(\Delta,P)$ and the  Riesz constants of the associated exponential systems. Additionally, this also links the types of exponential systems (Riesz basis, Riesz sequence, or frame) with the aspect ratio of the corresponding Fourier matrix (square, tall, or wide). There are several outstanding and natural questions which we plan to pursue in the future.
 
One direction is related to constructive perturbations and their explicit bounds. We present an observation about the principle minor of the DFT matrix and its implications to the stability of Riesz bases. 

\begin{Prop} \label{prop:instability}
	For any positive odd integer $N$, the basis $\{e^{2\pi i nx/N}\}_{n\in \Z}$ on $[0,N)$ has condition number 1, yet there is an explicit $\epsilon_N \colon \Z\to \R$ with $\|\epsilon_N \|_\infty=\frac1 2 - \frac 1{N+1}$ such that $\B_N=\{e^{2\pi i N^{-1}(n+\epsilon_N(n))x}\}_{n\in \Z}$ has condition number $\sqrt{N+1}$. 
\end{Prop}

The proof of this proposition is in \cref{proofs}. Notice that $\lim_{N\to\infty} \|\epsilon_N\|_\infty = \frac 1 2$  while the condition number  of $\B_N$ tends to infinity. This example also shows why in \cref{T-New-Kadec}, it is not possible to provide explicit frame bounds for the perturbed Riesz basis. 

This example naturally leads to the question of whether $\frac 12$ is the critical stability for Fourier matrices, or if $\frac 14$ is, like in the Kadec theorem. More specifically, we ask whether we can find a sequence of $N\in\N$ that goes to infinity and explicit perturbations $\epsilon_N$ of the frequencies of the $N\times N$ discrete Fourier matrix $[e^{2\pi i kj/N}]_{0\leq k\leq N-1\atop{0\leq j\leq N-1}}$ such that $\|\epsilon_N\|_\infty \to \frac 14$ and $\kappa_N\to\infty$, where $\kappa_N$ is the condition number of the perturbed matrix?  

A natural starting point is to examine a standard example which shows that $\frac 14$ is optimal for Kadec, given in \cite[Page 44]{Y} and originally proved in \cite{Levinson}. Let $\sign(\cdot)$ be the sign function with the convention that $\sign(0)=0$. Let $\epsilon\colon \Z\to\R$ be defined as $\epsilon(k)=-\frac 14\sign(k)$. When the frequencies of the usual exponential basis for $L^2([0,1])$ is perturbed by $\epsilon$, the resulting exponential system is no longer a basis. 

\begin{figure}[h]
	\includegraphics[width=0.5\textwidth]{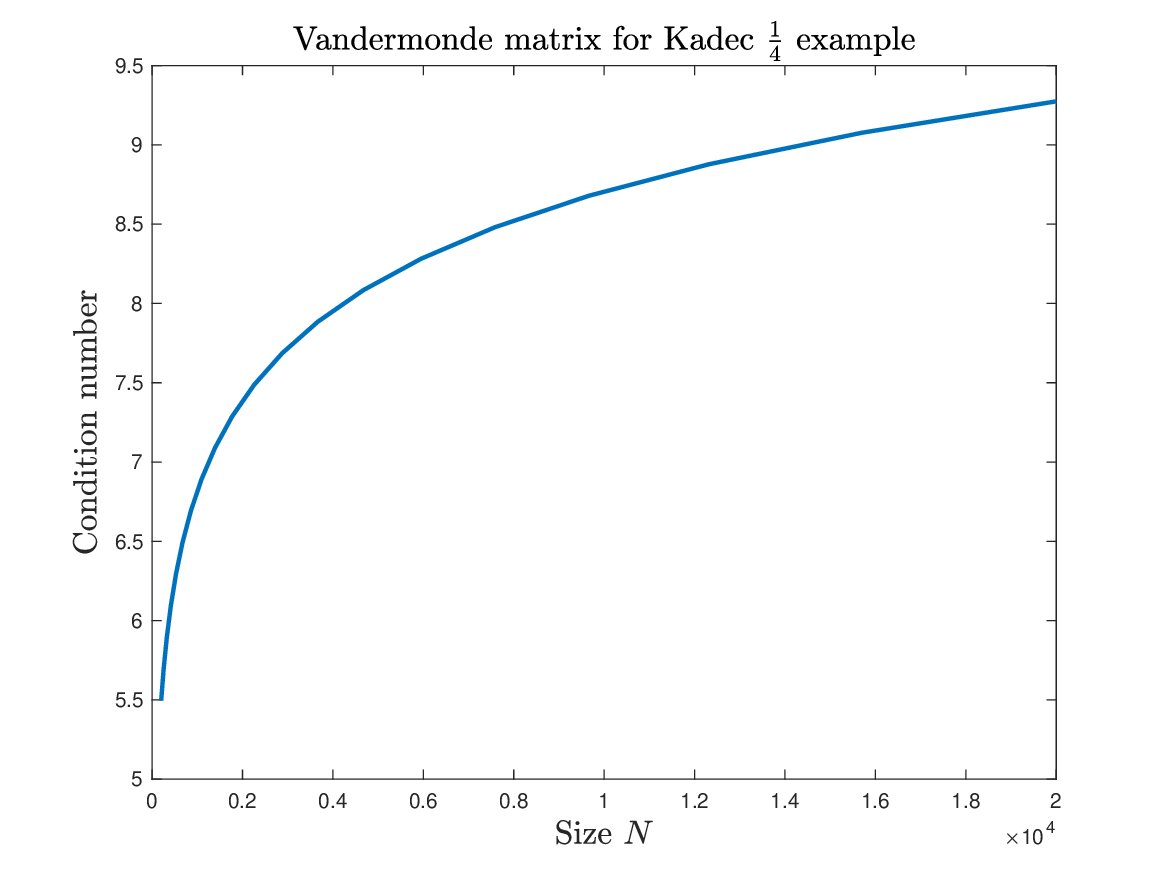}
	\caption{Condition number of $\bF_N'$ as a function of $N$.}
	\label{fig:Fcond}
\end{figure}

Unfortunately this example cannot be used for matrices; without explaining the full details, only periodic perturbations have a corresponding statement for Fourier matrices. Nonetheless, for a given odd integer $N=2m+1$, we can truncate $\epsilon$ and make it $N$-periodic on $\Z$ for arbitrary large $N$. Hence, let $\epsilon_N$ be the $N$-periodic sequence such that $\epsilon_N(k)=-\frac 14 \sign(k)$ for $k\in \{-m,\dots,m\}$. We examine the condition number of the perturbed Fourier matrix,
$$
\bF_N':=\Big[ \, e^{2\pi i j(k+\epsilon_N(k))/N } \, \Big]_{0\leq j,k<N}. 
$$ 
Figure \ref{fig:Fcond} displays its numerically computed condition number. The results of this experiment are inconclusive because while the condition number appears to grow in $N$, it is not clear whether it tends to infinity. The authors were only able to numerically compute the condition number up to matrices of size $20,000\times 20,000$ due to resource limitations. 

The stability  theorems that we have employed in this paper were originally proved for bases on rectangles in $\R^d$. We would like to derive more general and quantitative stability results  for exponential bases on other domains of $\R^d$. 
Specifically, we wonder whether Kadec-like stability  theorems could be proved for orthogonal exponential bases on   a domain $D\subset\R^n$  and what are the corresponding statements for matrices.  Recall that a bounded domain $D\subset\R^d$ is called {\it spectral} if $L^2(D)$ is  has an orthogonal exponential basis.  It was conjectured  by B. Fuglede \cite{F} that  $D$ is spectral if and only if it  tiles $\R^d$ by  translations.  The conjecture is false in both directions for general sets in  dimension $d\ge 2$, but it has been recently proved true for convex sets \cite{LM}.      

We also intend to explore the stability of sub-matrices of Fourier/Vandermonde matrices. They can naturally be associated to exponential systems on domains $T(P)$ as in \eqref{eq:TP}, where $T(P)$ is not necessarily a rectangle. To accomplish this goal, we aim to enhance our understanding of the stability of exponential bases on domains $T(P)$. Subsequently, we plan to refine \cref{prop:badstability2} and \cref{prop:badstability}. 

\section{Proof of \cref{T1-frame}} 
\label{sec:proofframe} 

\begin{proof}[Proof of \cref{T1-frame}] 
	{\it Proof of (a).}	Assume $L\ge N$ and that the   matrix $\G(\Delta,P)$ has rank $N$.
	We first derive a general relationship between $\B(\Delta)$ on $T(P)$ and matrix $\G(\Delta,P)$ that holds regardless of the relationship between $N$ and $L$. To simplify the resulting notation, let $\B:=\B(\Delta)$ and $T:=T(P)$.  Fix any $f\in L^2(T)$. For each $j\in \{1,\dots,L\}$, we obtain
	\begin{gather}\label{e3-frame-prod}
		\begin{split}
			\sum_{ \vec n\in\Z^d}  
			& \big| \l f, e^{ 2\pi i (\vec n+\vec \delta_j) \cdot x}  \r_{L^2(T)} \big|^2 
			= \sum_{ \vec n\in\Z^d} \Big| \sum_{k=1}^{N}\l f, e^{ 2\pi i (\vec n+\vec \delta_j) \cdot x}\r_{L^2( Q+\vec p_k)}\Big|^2
			\\ & 
			=  \sum_{\vec n\in\Z^d } \sum_{k,  k'=1}^{N}\l f, e^{ 2\pi i (\vec n+\vec \delta_j) \cdot x}\r_{L^2( Q+\vec p_k)}\overline{\l f, e^{ 2\pi i (\vec n+\vec \delta_j) \cdot x}\r_{L^2( Q+\vec p_{k'})}}.
		\end{split}
	\end{gather}
	Observe that
	\begin{gather} \label{e-Id1}
		\begin{split}
			&\l f,    e^{ 2\pi i (\vec n+\vec \delta_j) \cdot x}  \r_{L^2( Q +\vec p_k)}  =\int_{ Q +\vec p_k }  f(x) e^{ -2\pi i (\vec n+\vec \delta_j) \cdot x} \, dx \\ 
			&=	\int_{Q}  f(y+\vec p_k) e^{ -2\pi i (\vec n+\vec \delta_j)\cdot (y+\vec p_k)} \, dy 
			=
			e^{-2\pi i \vec \delta_j \cdot \vec p_k} \int_{Q } \tau_{\vec p_k} f(y) e^{ -2\pi i (\vec n+\vec \delta_j) \cdot y} \, dy \\    
			& = e^{-2\pi i \vec \delta_j \cdot \vec p_k} \l \tau_{\vec p_k}f , \ e^{  2\pi i (\vec n+\vec \delta_j)\cdot x}\r_{L^2(Q )},  
		\end{split}
	\end{gather}
	where we have defined $\tau_{\vec v}f(x)= f(x+\vec v)$. Combining \eqref{e3-frame-prod} and \eqref{e-Id1}, using that the set $\{e^{2\pi i (n+\vec \delta_j) \cdot x}\}_{\vec n\in\Z^d}$ is an orthonormal basis for $L^2(Q)$, and combining with the Parseval identity, we obtain  
	\begin{gather*}\label{e-Id2}
		\begin{split}
			\sum_{ \vec n\in\Z^d}  &    | \l f, e^{ 2\pi i (\vec n+\vec \delta_j) \cdot x}  \r_{L^2(T)}|^2  \\
			&=\sum_{k, k'=1}^{N} \!e^{ 2\pi i \vec \delta_j \cdot (\vec p_k-\vec p_{k'})}\!\!\sum_{n\in\Z^d }  \l \tau_{\vec p_k }f  ,   e^{  2\pi i (\vec n+\vec \delta_j)\cdot x}\r_{L^2(Q )}\overline{\l \tau_{\vec p_{k'}}f ,   e^{  2\pi i(\vec n+\vec \delta_j) \cdot x}\r_{L^2(Q )} } \\
			&=\sum_{k, k' =1}^{N }e^{ 2\pi i  \vec \delta_j \cdot (\vec p_k-\vec p_{k'})} \l  \tau_{\vec p_k} f, \, \tau_{\vec p_{k'}} f\r_{L^2(Q )} 
			=\sum_{k, k' =1}^{N } \beta_{k,k'} \l  \tau_{\vec p_kp} f,\     \tau_{\vec p_{k'}} f\r_{L^2(Q )},
		\end{split}  
	\end{gather*}
	where we have let
	\begin{equation*}\label{def-beta}
		\beta_{k,k'}
		= \sum_{j=1}^L e^{ 2\pi i \vec \delta_j \cdot (\vec p_{k'}-\vec p_k)},\quad 1\leq k,\, k'\leq N.
	\end{equation*}
	Let $\bB$ be the $N\times N$ matrix whose elements are the $\beta_{k,k'} $. A direct calculation shows that $\bB= \G^* \G$. Hence, for all $f\in L^2(T)$, it holds that 
	\begin{equation}\label{eq:Id3}
		\sum_{ \vec n\in\Z^d} | \l f, e^{ 2\pi i (\vec n+\vec \delta_j) \cdot x}  \r_{L^2(T)}|^2  
		=\sum_{k, k' =1}^{N } (\G^* \G)_{k,k'} \l  \tau_{\vec p_k} f,\     \tau_{\vec p_{k'}} f\r_{L^2(Q)}. 
	\end{equation}
	
	Now we concentrate on the upper bound in part (a). Using the variation characterization of the maximum eigenvalues of a Hermitian matrix in \eqref{eq:Id3}, we see that     
	\begin{gather} \label{eq:Id4}
		\begin{split}
			\sum_{\vec n\in\Z^d} | \l f, e^{ 2\pi i (\vec n+\vec \delta_j) \cdot x}  \r_{L^2(T)}|^2
			&= \int_Q \sum_{k, k' =1}^{N } (\G^* \G)_{k,k'} \tau_{\vec p_k} f(x) \tau_{\vec p_{k'}} f(x) \, dx \\
			&\leq \sigma_1^2(\G) \int_{Q }\sum_{k=1}^{N }  |\tau_{\vec p_k} f(x)|^2 \, dx 
			= \sigma_1^2(\G) \|f\|_{L^2(T)}^2.
		\end{split}
	\end{gather}
	This establishes that $\sigma_1^2(\G)$ is a valid upper frame bound for $\B$ on $T$. To show this is optimal, let $\phi\in \C^N$ be a unit norm vector such that $\bB\phi = \lambda_1(\bB) \phi$. Let $f\in L^2(T)$ such that $\tau_{\vec p_k} f(x)=\phi_k$ for all $x\in Q$. Doing so shows that inequality \eqref{eq:Id4} is attained by $f$. 
	
	Now we concentrate on the lower bound in part (a). We repeat a similar argument as in \eqref{eq:Id4} to see that
	$$
	\sum_{\vec n\in\Z^d} | \l f, e^{ 2\pi i (\vec n+\vec \delta_j) \cdot x}  \r_{L^2(T)}|^2
	\geq \sigma_N^2(\G) \int_{Q }\sum_{k=1}^{N }  |\tau_{\vec p_k} f(x)|^2 \, dx 
	= \sigma_N^2(\G) \|f\|_{L^2(T)}^2.
	$$	
	This proves that $\sigma_N^2(\G)$ is a valid lower frame constant for $\B$ on $T$. To prove that this is optimal, let $\psi\in \C^N$ be a unit norm vector such that $\bB\phi = \lambda_N(\bB) \phi$ and $\bB\psi=\lambda_N(\bB)\psi$. Let $g\in L^2(T)$ such that $\tau_{\vec p_k} g(x)=\psi_k$ for all $x\in Q$ and $k\in \{1,\dots,N\}$. Doing so shows that the above lower bound is attained by $g\in L^2(T)$. 
	
	\medskip \noindent
	{\it Proof of part (b)}.  We now consider the case $L\leq N$; assume  for simplicity that $d=1$, but the proof extends to higher dimensions (see \cite{DS} for additional details). 
	
	For the sake of uniformity of notation, we will denote with $Q$ the unit segment $[-\frac 12, \frac 12)$.    We assume that $\G$ has rank $L$ and we  show that $\B(\Delta)$ is a Riesz sequence in $L^2(T)$  by proving   that there exist $A, B>0$ for which  
	\begin{equation}\label{e2-cond-lin-ind} 
		A  \sum_{j=1}^{L}\sum_{  n\in\Z } | a_{j,  n }|^2 \leq  
		\Big\| \sum_{j=1}^{L}\sum_{  n\in\Z } a_{ j,  n } \, e^{ 2\pi i  (n+\delta_j)x}  \Big\|_{L^2(T)}^2 
		\leq  B  \sum_{j=1}^{L}\sum_{  n\in\Z } | a_{  n,j}|^2 \end{equation} 
	for every  set of coefficients $\{a_{j,  n}\}_{j\leq L\atop{n\in \Z}}\subset\C$ which are non-zero for only finitely many pairs $(j,n)$. To this aim, we use a family of  generalized  Hilbert transform operators developed in \cite{DS}.  For a given $t\in\R$, we let $H_t:\ell^2(\Z)\to\ell^2(\Z)$
	$$ 
	(H_t(a))_m= 
	\begin{cases}  \, \dsize \frac{ \sin ( \pi t)}{ \pi}\sum_{  n\in\Z }      \frac{a_n  }{m-   n  +t} 
	&\text{if } t\not\in\Z, \\
	\ (-1)^t a_{m+t} &\text{if } t \in\Z. \end{cases}
	$$ 
	In \cite{DS} it is  proved that the family $\{H_t\}_{t\in\R}$   is a strongly continuous  Abelian group   of isometries  in $ \ell^2(\Z)$; in particular, it is   proved that $   H_{s}\circ H_t = H_{s+t} $, that $H_s^{-1}= H_s^{*}=H_{-s}$  and $\|H_s(a)\|_{\ell^2(\Z)}= \|a\|_{\ell^2(\Z)}$  for every  $a\in \ell^2(\Z) $. Lemma 4.1 in \cite{DC} shows that for every $a, b\in\ell^2(\Z )$ and every $s,t\in \R $ and $v\in\Z$, we have  
	\begin{equation}\label{e2-inner-Sj} 
		\Big\l \sum_{  n\in\Z } a_{  n} e^{2\pi i\, (n+  s)x},\ \sum_{ m\in\Z } b_{ m}e^{2\pi i(m +  t)x} \Big\r_{L^2(Q+v)}
		=  e^{2\pi i (s -  t)v  }\l  H_{  t}( \alpha ), H_{  s}( \beta)\r_{\ell^2(\Z )}, 
	\end{equation}
	where $\alpha  = ((-1)^{n }   a_{  n})_{  n\in\Z }$ and $ \beta = ((-1)^{m } b_{  m})_{  m\in\Z }$. 
	
	Let 
	$  \dsize S_j  =\sum_{  n \in\Z } a_{j,  {n }} e^{2\pi i\, (n+\delta_j)x  }  .$ 
	In view of \eqref{e2-inner-Sj}, the  central sum in \eqref{e2-cond-lin-ind} equals 
	\begin{align*}
		\|S_1+\cdots+S_L\|^2_{L^2(T)}  
		&=   \sum_{k=1}^N\sum_{1\leq i, j\leq L}\langle S_i, S_j\rangle_{L^2(Q+p_k)} \\
		&=  \sum_{k=1}^N\sum_{1\leq i, j\leq L} e^{2\pi i ( \delta_i-\delta_j)p_k} \l H_{ - \delta_i }(a_i),  H_{-\delta_j}(a_j)\r_{\ell^2(\Z)} .
	\end{align*}
	Letting $b_s:=H_{-\delta_s }(a_s)$ and observing that  $\sum_{k=1}^L e^{2\pi i (\delta_i-  \delta_j)p_k}= (\G \G^*)_{i,j}$, we can write 
	\begin{align*}
		\|S_1+\cdots+S_L\|^2_{L^2(T)}
		%&=\sum_{k=1}^N\sum_{i,j=1}^{L }  e^{2\pi i (\delta_i-\delta_j)p_k} \l H_{-\delta_i }( a_i),  H_{-\delta_j}(a_j)\r_{\ell^2(\Z )} \\
		&=\sum_{i,j=1}^{L } (\G \G^*)_{i,j} \l  b_i,\,  b_j\r
		=  \sum_{n\in\Z} \sum_{i,j=1}^{L } (\G \G^*)_{i,j} \,  b_{i,n} \, \overline b_{j,n}.
	\end{align*} 
	By the variational characterization of the maximum/minimum  eigenvalues of a Hermitian matrix, we obtain 
	$$
	\sigma_L^2(\G)  \sum_{j=1}^L \|  b_j\|^2_{\ell^2(\Z )}
	\leq \sum_{n\in\Z}\sum_{i,j=1}^{L} (\G \G^*)_{i,j} \, b_{i,n} \, \overline b_{j,  n}
	\leq \sigma_1^2(\G) \sum_{j=1}^L \| b_j\|^2_{\ell^2(\Z )}. 
	$$
	The central sum in the expression above is $\|S_1+\cdots +S_{L}\|^2_{L^2(Q)}$.
	Since   $ b_j=H_{- \delta_j}(a_j)$ and   the $H_{ t}$   are invertible isometries,  we have   $\|  b_j\|_{\ell^2(\Z )}=  \| a_j\|_{\ell^2(\Z )}$,  and thus,
	$$
	\sigma_L^2(\G)  \sum_{j=1}^L \|a_j\|^2_{\ell^2(\Z)} 
	\leq \|S_1+ \cdots +S_{L}\|^2_{L^2(Q)} 
	\leq  \sigma_1^2(\G)  \sum_{j=1}^L  \|a_j\|^2_{\ell^2(\Z)}. 
	$$
	This proves \eqref{e2-cond-lin-ind} with $A= \sigma_L^2(\G)$ and $B=\sigma_1^2(\G)$ as valid Riesz constants. The proof of optimality is similar to that in part (a) and exploits that $H_t$ is an isometry. 
\end{proof}

\section{Other proofs} 
\label{proofs}

\begin{proof}[Proof of \cref{T1-special-delta}]
	%We can assume $K=1$ without loss of generality. 
	 We need to prove that  the sub-matrix  $\G_N$  of $\G$ formed with the first $N$ rows of $\G$ is nonsingular if and only of  \eqref{e2-cond=delta} holds. \cref{T1-frame} yields that $\B$ is a frame. The matrix $\G_N$ is   Vandermonde and its determinant is 
	$$\mbox{det}\G_N=\prod_{1\leq k <k'\leq N} (e^{2\pi i \l \vec p_k,\, \vec\delta\r}-e^{2\pi i \l \vec p_{k'},\, \vec\delta\r});
	$$
	Thus, $\G_N$ is nonsingular if and only if $e^{2\pi i\l \vec p_k-\vec p_{k'},\, \vec\delta\r}\ne 1$, or $\l \vec p_k-\vec p_{k'},\, \vec\delta\r\not\in\Z$.
\end{proof}

\medskip
\begin{proof}[Proof of Theorem \ref{T-New-Kadec}]	
	For each $s=1,\dots, d$, we let  $\delta_{j_k}:= \frac{j_k +\epsilon_k(j_k)}{M_k}$.  The set, $$
	\B_k=\bigcup_{0\leq j_k\leq M_k-1} \big\{e^{ 2\pi i (n_k+ \delta_{j_k} )x_k } \big\}_{  n \in\Z },$$ 
	is a perturbation of the standard exponential basis for $L^2([0, M_k))$. By \cite[Corollary 5.1]{DC}, we have $\B_k$ is a basis if and only if $\delta_{j_k}-\delta_{i_k }\not\in\Z$, i.e., 
	if and only if  ${M_k}^{-1} (\epsilon_k(i_k)-\epsilon_k(j_k) +i_k-j_k)\not\in\Z$ whenever $i_k\ne j_k$. The exponential system $\B$ is the tensor product of $\B_1,\dots,\B_d$. Thus, $\B$ is an exponential basis for  $L^2(R(\vec M))$ if and only if each $\B_k$  is a basis for $L^2([0, M_k))$  and   by \cref{T1-special-delta}, if and only if condition \eqref{e-cond1} holds.
\end{proof}

\begin{proof}[Proof of \cref{prop:badstability2}]
	Let $\bA$ and $\bE$ be matrices of size $L\times N$. By  the triangle inequality, 
	$
	\sigma_1(\bA+\bE)
	\leq \sigma_1(\bA)+\sigma_1(\bE). 
	$
	On the other hand, the Weyl's inequality (see e.g. \cite[page 474]{HJ})  yields	$
	\sigma_N(\bA+\bE)\geq \sigma_N(\bA)-\sigma_1(\bE).$
	
	We let $\bA=\bF(\Omega,U)$ and $\bE=\bF(\Omega',U)-\bF(\Omega,U)$, respectively. It remains to upper bound $\sigma_1(\bE)$. Letting $\|\bE\|_F $ be the Frobenius norm,  we have $\sigma_1(\bE)\leq \|\bE\|_F$. By the hypotheses, for each $j\in  \{1,\dots,L\}$ and $\vec u \in U$, we have
	$$
	|e^{2\pi i \vec \omega_j \cdot \vec u}-e^{2\pi i \vec \omega_j' \cdot \vec u}|
	=|1-e^{2\pi i (\vec \omega_j'-\vec \omega_j) \cdot \vec u}|
	\leq 2\pi |(\vec \omega_j'-\vec \omega_j) \cdot \vec u|
	\leq 2\pi \epsilon \|u\|_{p'}
	\leq \pi \epsilon d^{1/p'}.  
	$$
	We have  used the mean value theorem and that $\|\vec u\|_2\leq d^{1/p'} \|u\|_\infty \leq \frac 12 d^{1/p'}$. This implies 
	$$
	\|\bE\|_F
	\leq \left(\sum_{\vec u\in U} \sum_{j=1}^L |e^{2\pi i \vec \omega_j \cdot \vec u}-e^{2\pi i \vec \omega_j' \cdot \vec u}|^2 \right)^{1/2}
	\leq \pi d^{1/p'} \sqrt{L |U|} \epsilon. 
	$$	
\end{proof}  

\begin{proof}[Proof of \cref{prop:badstability}]
	Following the same strategy as in the proof of \cref{prop:badstability2}, we let $\bA=\bF(\Omega,U)$ and $\bE=\bF(\Omega,V)-\bF(\Omega,U)$, and we   estimate $\sigma_1(\bE)$. By the mean value theorem, for each $\vec \omega\in\Omega$ and $k\in \{1,\dots,N\}$, we have 
	$$
	|e^{2\pi i \vec \omega \cdot \vec u_k'} - e^{2\pi i \vec \omega \cdot \vec u_k}|
	=|1-e^{2\pi i \vec \omega \cdot (\vec u_k-\vec u_k')}|
	\leq 2\pi |\vec \omega\cdot (\vec u_k-\vec u_k')|
	\leq 2\pi \epsilon \|\vec \omega\|_{p'}. 
	$$
	This implies 
	$$
	\sigma_1(\bE)
	\leq \|\bE\|_F
	\leq \left(\sum_{k=1}^N \sum_{\omega\in\Omega} |e^{2\pi i \vec \omega \cdot \vec u_k'} - e^{2\pi i \vec \omega \cdot \vec u_k}|^2\right)^{1/2}
	%\leq 2\pi \sqrt N \epsilon \left(\sum_{\omega\in \Omega} \|\omega\|_{p'}^2\right)^{1/2}
	\leq C_{\Omega,p} \sqrt N \epsilon. 
	$$	
\end{proof}

\begin{proof}[Proof of \cref{prop:instability}]
	Let $N=2m+1$. As discussed earlier, the associated $\G$ matrix to $\B:=\{e^{2\pi i nx/N}\}_{n\in \Z}$ on $[0,N)$ is simply the DFT matrix of size $N\times N$. Let $\epsilon$ be the sequence  $\{-\frac{j}{N+1}\}_{j=-m,\dots,m}$ extended to be $N$ periodic on $\Z$ and notice that 
	$
	\|\epsilon\|_\infty =\frac{m}{N+1} 
	= \frac1 2 - \frac 1{N+1}.
	$
	Next, we perturb $\B$ to the sequence
	$$
	\B'
	=\Big\{e^{2\pi i N^{-1}(n+\epsilon(n))x} \Big\}_{n\in \Z}
	=\bigcup_{j=-m}^{m} \Big\{e^{2\pi i (n+\frac{j}{N+1})x} \Big\}_{n\in \Z}. 
	$$
	Notice that $\B'=\B(\Delta')$, where $\Delta'=\{\frac j{N+1}\}_{j=0,\dots,N-1}$. Thus the matrix associated to the sequence $\B'$ on the interval $[0,N)$ is
	$$
	\G(\Delta',P)
	=\Big[ e^{2\pi i jk/(N+1)}\Big]_{j=-m,\dots,m \,\atop{ k=0,\dots,N-1}}. 
	$$
	After re-indexing, we see that $\G(\Delta',P)$ coincides with the leading $N\times N$ submatrix $\bA$ of the $(N+1)\times (N+1)$ DFT matrix $\bF$. 
	 
	We claim that the singular values of $\bA$ are precisely $1$ and $\sqrt{N+1}$, with multiplicities $1$ and $N-1$ respectively. To see this, let $\bU$ and $b^*$ be the first $N$ columns and last row of $\bF$, respectively. Since $\bU$ consists of orthogonal columns each with $\ell^2$ norm $\sqrt{N+1}$, we have  
	$
	(N+1) \bI_{N+1} = \bU^* \bU = \bA^*\bA + bb^*.
	$
	This equation implies that $b$ is an eigenvector of $\bA^*\bA$ with eigenvalue $1$ and any $v$ orthogonal to $b$ is also an eigenvector of $\bA^*\bA$ with eigenvalue $N+1$. By selecting an orthonormal basis for the orthogonal complement of $b$, this establishes the claim. Consequently, the condition number of $\bA'$, $\G(\Delta',P)$, and $\B'$ are each $\sqrt{N+1}$. 
\end{proof}

\section*{Acknowledgments} 

The authors  wish to thank Azita Mayeli for   her  helpful insights during an early stage of this project.

WL is partially supported by NSF-DMS Award \#2309602 and a PSC--CUNY grant.

\end{document}